\documentclass[
	english
]{scrartcl}

\usepackage[T1]{fontenc}
\usepackage[utf8]{inputenc}

\usepackage{amsmath,amsfonts,amsthm,amssymb}
\usepackage[colorlinks,citecolor=blue]{hyperref}
\usepackage{doi,natbib}
\usepackage{natbib}
\usepackage{cancel}
\usepackage{relsize}

\usepackage{changes}

\usepackage{enumitem}
\setlist[enumerate]{label=(\alph*)}

\usepackage[figure]{hypcap}
\usepackage{graphicx}
\graphicspath{{./img/}}
\usepackage{subcaption} 

\usepackage{preamble}

\usepackage{marginnote}

\usepackage[capitalise,nameinlink]{cleveref}
\allowdisplaybreaks

\numberwithin{equation}{section}

\newcommand\norm[1]{\left\Vert#1\right\Vert}

\newcommand\innerprod[2]{\left\langle #1, #2\right\rangle}
\newcommand\ninnerprod[2]{\langle #1, #2\rangle}

\newcommand\N{\mathbb{N}}
\newcommand\R{\mathbb{R}}

\newcommand\B{\mathbb{B}}

\newcommand{\dist}{\operatorname{dist}}

\DeclareMathAlphabet{\mathpzc}{OT1}{pzc}{m}{it}

\newtheorem{theorem}{Theorem}[section]
\newtheorem{lemma}[theorem]{Lemma}

\newtheorem{corollary}[theorem]{Corollary}

\definecolor{mygreen}{rgb}{0.0,0.7,0.0}
\definecolor{mybrown}{rgb}{0.5,0.5,0.0}

\begin{document}

\title{A simple proof of second-order sufficient optimality conditions in nonlinear semidefinite optimization}
\author{%
	Patrick Mehlitz%
	\footnote{%
		Brandenburgische Technische Universit\"at Cottbus--Senftenberg,
		Institute of Mathematics,
		03046 Cottbus,
		Germany,
		\email{mehlitz@b-tu.de},
		\url{https://www.b-tu.de/fg-optimale-steuerung/team/dr-patrick-mehlitz},
		\orcid{0000-0002-9355-850X}%
		}
	}

\publishers{}
\maketitle

\begin{abstract}
	In this note, we present an elementary proof for a well-known 
	second-order sufficient optimality
	condition in nonlinear semidefinite optimization which does not rely on the enhanced
	theory of second-order tangents.
	Our approach builds on an explicit elementary
	computation of the so-called second subderivative
	of the indicator function associated with the semidefinite cone which recovers the
	best curvature term known in the literature.
\end{abstract}


\begin{msc}	
	\mscLink{49J52}, \mscLink{49J53}, \mscLink{90C22}
\end{msc}

\section{Introduction}\label{sec:introduction}

Second-order sufficient optimality conditions play a significant role in the theory of
nonlinear optimization. Among others, their validity guarantees stability of the
underlying strict local minimizer with respect to perturbations of the data, and this opens a 
way in order to show local fast convergence of diverse types of numerical solution
algorithms, including augmented Lagrangian, sequential quadratic programming, and
Newton-type methods.

Geometric constraints of type
\begin{equation}\label{eq:geometric_constraints}
	F(x)\in C,
\end{equation}
where $F\colon\mathbb X\to\mathbb Y$ is a twice continuously differentiable mapping between
Euclidean spaces $\mathbb X$ and $\mathbb Y$, and $C\subset\mathbb X$ is a closed, convex set,
provide a rather general paradigm 
for the modeling of diverse popular constraint systems in nonlinear optimization. 
It has been well-recognized in the past that second-order optimality conditions in
constrained optimization depend on the second derivative of the objective function as well as
the curvature of the feasible set. In the presence of constraints of type
\eqref{eq:geometric_constraints}, the latter can be described in terms of the second
derivative of $F$ and the curvature of $C$. Thus, associated second-order optimality
conditions do not only comprise the second derivative of a suitable Lagrangian function,
but a so-called curvature term associated with $C$ pops up as well.
In case where $C$ is a polyhedral set, this curvature term vanishes, and one obtains 
very simple second-order conditions as they are known from standard nonlinear programming,
see \cite{BenTal1980,McCormick1967}. In more general cases, however, a suitable tool to
keep track of the curvature of $C$ has to be used to formulate a suitable curvature term.
Classically, the support function of a (local) second-order tangent approximation of $C$
has been exploited for that purpose, see \cite{BonnansCominettiShapiro1999,BonSh00}, and this
exemplary led to second-order optimality conditions in nonlinear second-order cone and
semidefinite optimization, see \cite{BonnansRamirez2005,Shapiro1997}. 
However, we would like to mention here that the proofs in these papers are far from being
elementary since the calculus of second-order tangents is a rather challenging task.
With the aid of a generalized notion of support functions, 
the approach via second-order tangents can be further generalized to situations 
where $C$ is not convex anymore, see \cite{GfrererYeZhou2022}.
Another less popular approach to curvature terms has been promoted recently in
\cite{BeGfrYeZhouZhang1,MohammadiMordukhovichSarabi2021,ThinhChuongAnh2021} where the 
so-called second subderivative, see \cite{Rockafellar1989}, of the indicator function
of $C$ has been used for that purpose. 
This tool yields promising results even in infinite-dimensional spaces,
see \cite{ChristofWachsmuth2018,WachsmuthWachsmuth2022}. 
The approach via second subderivatives is particularly suitable for the
derivation of second-order sufficient optimality conditions due to the underlying
calculus properties of second subderivatives, see \cite{BenkoMehlitz2022c} for a recent study.
Second-order sufficient conditions obtained from this approach have been shown to serve
as suitable tools for the local convergence analysis of solution algorithms associated with
challenging optimization problems based on variational analysis, 
see \cite{HangMordukhovichSarabi2022,HangSarabi2021,Sarabi2022}.
In this note, we aim to popularize the approach using second subderivatives even more
by presenting an application in nonlinear semidefinite optimization.

Thus, let us focus on the special situation where $\mathbb Y:=\mathbb S^m$ equals the
space of all real symmetric $m\times m$-matrices and $C:=\mathbb S^m_+$ is the cone of all 
positive semidefinite matrices.
The tightest second-order sufficient condition in nonlinear semidefinite optimization
we are aware of has been established by Shapiro and 
can be found in \cite[Theorem~9]{Shapiro1997}. Its proof heavily 
relies on technical arguments which exploit second-order directional 
differentiability of the smallest eigenvalue of a positive semidefinite matrix and
calculus rules for second-order tangent sets.
Later, several authors tried to recover or enhance this result using reformulations
of the original problem.
In \cite{Forsgren2000}, the author obtained a related second-order sufficient condition based
on a localized Lagrangian and some technical arguments via Schur's complement.
The authors of \cite{LourencoFukudaFukushima2018} applied the squared slack variable technique
to semidefinite optimization problems and obtained second-order sufficient conditions in the
presence of so-called strict complementarity.
In \cite{Jarre2012}, strict complementarity and a second-order
constraint qualification are needed to recover Shapiro's original second-order sufficient
condition based on a simplified technique. 
Further results about second-order optimality conditions in nonlinear semidefinite
optimization such as a strong second-order sufficient condition and a weak second-order
necessary condition can be found in \cite{FukudaHaeserMito2020, Sun2006}.
The validation of second-order sufficient conditions in the papers 
\cite{Forsgren2000,Jarre2012,LourencoFukudaFukushima2018}
is much simpler than the strategy used in \cite{Shapiro1997}.
However, these approaches either do not recover the original result 
from \cite{Shapiro1997} in full generality, i.e., additional conditions are
postulated to proceed, or the analysis still makes some technical preliminary
considerations necessary.
Here, we simply compute the second subderivative of the indicator function associated
with the positive semidefinite cone in order to recover the result from \cite{Shapiro1997}
in elementary way.
	Let us note that this calculation already has been done in 
	\cite[Example~3.7]{MohamadiSarabi2020}, but the arguments presented there
	are not self-contained and exploit involved variational properties of
	eigenvalue functions, see \cite{Torki1999}.
	In contrast, our calculations are completely elementary.

The remainder of this note is structured as follows.
In \cref{sec:preliminaries}, we summarize the notation used in this paper and recall the
definitions of some variational tools which we are going to exploit.
We present an abstract second-order sufficient optimality condition for nonlinear
semidefinite optimization problems in \cref{sec:SOSC} which comprises the second subderivative
of the indicator function of the semidefinite cone 
as the curvature term and can be distilled from a much more general result recently proven in
\cite{BeGfrYeZhouZhang1,BenkoMehlitz2022c}. Then, by explicit computation of the appearing
second subderivative, we specify this result in terms of initial problem data
and recover the results from \cite{BonSh00,Shapiro1997}. 
Some concluding remarks close the paper in \cref{sec:conclusions}.

\section{Preliminaries}\label{sec:preliminaries}

The notation used in this note is fairly standard and follows
\cite{BonSh00,RoWe98}.

\subsection{Basic notation}

By $\R^n_+$, we denote the nonnegative orthant of $\R^n$.
Let $\R^{m\times n}$ be the set of all rectangular matrices with
$m$ rows amd $n$ columns, and $O$ the all-zero matrix of appropriate dimensions. 
An Euclidean space $\mathbb X$, i.e., a finite-dimensional Hilbert space,
will be equipped with the inner product 
$\innerprod{\cdot}{\cdot}\colon\mathbb X\times\mathbb X\to\R$ and the
associated induced norm $\norm{\cdot}\colon\mathbb X\to[0,\infty)$.
For arbitrary $\bar x\in\mathbb X$ and $\varepsilon>0$, 
$\mathbb B_\varepsilon(\bar x):=\{x\in\mathbb X\,|\,\norm{x-\bar x}\leq\varepsilon\}$
represents the closed $\varepsilon$-ball around $\bar x$.
The space of all real symmetric $n\times n$-matrices $\mathbb S^n$ 
is equipped with the Frobenius inner product given by
\[
	\forall A,B\in\mathbb S^n\colon\quad
	\innerprod{A}{B}
	:=
	\operatorname{trace}(AB)
	=
	\sum_{i=1}^n\sum_{j=1}^nA_{ij}B_{ij}
\]
and the associated induced Frobenius norm.

For an arbitrary Euclidean space $\mathbb X$ and some nonempty, convex set $A\subset\mathbb X$,
we use
\begin{align*}
	A^\circ&:=\{x^*\in\mathbb X\,|\,\forall x\in A\colon\,\ninnerprod{x^*}{x}\leq 0\},\\
	A^\perp&:=\{x^*\in\mathbb X\,|\,\forall x\in A\colon\,\ninnerprod{x^*}{x}=0\}
\end{align*}
in order to denote the polar cone of $A$, which is always a closed, convex cone,
and the annihilator of $A$, which is a subspace of $\mathbb X$.
The distance function $\dist_A\colon\mathbb X\to\R$ of $A$ is given by
\[
	\forall x\in\mathbb X\colon\quad
	\dist_A(x):=\inf\{\norm{y-x}\,|\,y\in A\}.
\]
For $\bar x\in A$, we make use of
\begin{align*}
	\mathcal T_A(\bar x)
	:=
	\left\{
		u\in\mathbb X\,\middle|\,
			\begin{aligned}
				&\exists\{t_k\}_{k\in\N}\subset(0,\infty)\,
					\exists\{u_k\}_{k\in\N}\subset\mathbb X\colon\\
				&\quad t_k\downarrow 0,\,u_k\to u,\,\bar x+t_ku_k\in A\,\forall k\in\N
			\end{aligned}
	\right\}
\end{align*}
in order to represent the tangent (or Bouligand) cone to $A$ at $\bar x$.
The associated polar cone, i.e.,
\[
	\mathcal N_A(\bar x):=\mathcal T_A(\bar x)^\circ
\]
is the normal cone to $A$ at $\bar x$. 
Note that $\mathcal T_A(\bar x)$ and $\mathcal N_A(\bar x)$ are closed, convex cones.

For a twice continuously differentiable mapping $F\colon\mathbb X\to\mathbb Y$ between
Euclidean spaces $\mathbb X$ and $\mathbb Y$ as well as some point $\bar x\in\mathbb X$,
$F'(\bar x)\colon\mathbb X\to\mathbb Y$ is the linear operator which represents
the first derivative of $F$ at $\bar x$.
Similarly, $F''(\bar x)\colon\mathbb X\times\mathbb X\to\mathbb Y$ is the bilinear
mapping which represents the second derivative of $F$ at $\bar x$.
Partial derivatives are denoted in analogous way.

Finally, for a lower semicontinuous function $\varphi\colon\mathbb X\to\R\cup\{\infty\}$,
some $\bar x\in\mathbb X$ such that $\varphi(\bar x)<\infty$, and some
$x^*\in\mathbb X$, the function 
$\mathrm d^2\varphi(\bar x,x^*)\colon\mathbb X\to\R\cup\{-\infty,\infty\}$ given by
\[
	\forall u\in\mathbb X\colon\quad
	\mathrm d^2\varphi(\bar x,x^*)(u)
	:=
	\liminf\limits_{t\downarrow 0,\,u'\to u}
	\frac{\varphi(\bar x+tu')-\varphi(\bar x)-t\ninnerprod{x^*}{u'}}{t^2/2}
\]
is referred to as the second subderivative of $\varphi$ at $\bar x$ with $x^*$.
The recent study \cite{BenkoMehlitz2022c} reports on the calculus of this variational
tool and its usefulness for the derivation of second-order optimality conditions in
nonlinear optimization, and these findings can be partially extended even to
infinite-dimensional situations, see \cite{ChristofWachsmuth2018,WachsmuthWachsmuth2022}. 
Here, we are particularly interested in the second subderivative
of indicator functions $\delta_A\colon\mathbb X\to\R\cup\{\infty\}$, associated with
closed, convex sets $A\subset\mathbb X$, given by
\[
	\forall x\in\mathbb X\colon\quad
	\delta_A(x)
	:=
	\begin{cases}
		0	&	x\in A,\\
		\infty	&	x\notin A.
	\end{cases}
\]
For this particular function, the definition of the second subderivative yields
\[
	\forall u\in\mathcal T_A(\bar x)\colon\quad
	\mathrm d^2\delta_A(\bar x,x^*)(u)
	=
	\liminf\limits_{\substack{t\downarrow 0,\,u'\to u\\\bar x+tu'\in A}}
	-\frac{2\ninnerprod{x^*}{u'}}{t},
\]
and one can easily check that $\mathrm d^2\delta_A(\bar x,x^*)(u)=\infty$ if
$u\notin\mathcal T_A(\bar x)$ or $\innerprod{x^*}{u}<0$.
In case where $u\in\mathcal T_A(\bar x)$ and $\ninnerprod{x^*}{u}>0$, 
$\mathrm d^2\delta_A(\bar x,x^*)(u)=-\infty$ holds.
Thus, only the case $u\in\mathcal T_A(\bar x)\cap\{x^*\}^\perp$ is interesting.
In turn, for given $\bar x\in A$ and $u\in\mathcal T_A(\bar x)$, the consideration
of the second subderivative is only reasonable if $x^*\in\mathcal N_A(\bar x)\cap\{u\}^\perp$.

\subsection{Matrix analysis}\label{sec:matrix_analysis}

In order to carry out our analysis related to the cone of all positive semidefinite
matrices, we need to introduce some further notation first.
Fix some $m\in\N$ such that $m\geq 2$.
By $\mathbb S^m_+,\mathbb S^m_-\subset\mathbb S^m$, we denote the cones of
all positive semidefinite and negative semidefinite matrices, respectively.
For each matrix $Y\in\mathbb S^m_+$, there exists an orthogonal matrix
$P\in\R^{m\times m}$ such that $Y=P^\top MP$ where $M\in\R^{m\times m}$ is
the diagonal matrix whose diagonal is made of the eigenvalues of $Y$, 
ordered non-increasingly. We refer to this representation as an 
ordered eigenvalue decomposition of $Y$. Throughout the paper, we
will denote the index sets of (row) indices of $M$ associated with the positive and zero
eigenvalues of $Y$ by $\pi$ and $\omega$, respectively.
For later use, let us also mention that $Y^\dagger =P^\top M^\dagger P$ holds for
the Moore--Penrose pseudoinverse of $Y$, 
and that $M^\dagger$ results from $M$ by inverting its positive diagonal elements.
For arbitrary matrices $A\in\mathbb S^m$ and index sets $I,J\subset\{1,\ldots,m\}$, 
we use $A_{IJ}$ to denote the matrix which results from $A$ by deleting those rows and
columns whose indices do not belong to $I$ and $J$, respectively.
Furthermore, we set $A^P:=PAP^\top$ and $A^P_{IJ}:=(A^P)_{IJ}$.

In \cite[Section~5.3.1]{BonSh00}, the formula
\begin{equation}\label{eq:tangent_cone_semidefinite}
	\mathcal T_{\mathbb S^m_+}(Y)
	=
	\left\{V\in\mathbb S^m\,\middle|\,V^P_{\omega\omega}\in\mathbb S^{|\omega|}_+\right\}
\end{equation}
has been established.
Furthermore, \cite[Section~4.2.4]{HiriartUrrutyMalick2012} gives
\begin{equation}\label{eq:normal_cone_semidefinite}
	\mathcal N_{\mathbb S^m_+}(Y)
	=
	\left\{Y^*\in\mathbb S^m\,\middle|\,(Y^*)^P_{\pi\pi}=O,\,(Y^*)^P_{\pi\omega}=O,\,
		(Y^*)^P_{\omega\omega}\in\mathbb S^{|\omega|}_-\right\}.
\end{equation}

	In the course of this note, we will need a criterion for semidefiniteness of block matrices.
	The following lemma is taken from \cite[Appendix~A.5.5]{BoydVandenberghe2004}.
	\begin{lemma}\label{lem:definiteness_block_matrix}
		Let $m_1,m_2\in\N$ be positive integers.
		Furthermore, let $A\in\mathbb S^{m_1}_+$ be positive definite, and let
		$B\in\R^{m_1\times m_2}$ as well as $C\in\mathbb S^{m_2}$ be arbitrarily chosen.
		For $m:=m_1+m_2$, we consider the block matrix
		\[
			M:=\begin{bmatrix}A&B\\B^\top&C\end{bmatrix}\in\mathbb S^m.
		\]
		Then $M\in\mathbb S^m_+$ is equivalent to $C-B^\top A^{-1}B\in\mathbb S^{m_2}_+$.
	\end{lemma}

\section{Second-order sufficient optimality conditions in nonlinear semidefinite optimization}\label{sec:SOSC}

Let $m\in\N$ such that $m\geq 2$ be fixed.
Throughout the section, we consider the nonlinear semidefinite optimization problem
\begin{equation}\label{eq:optimization_problem}\tag{NSDP}
	\min\{f(x)\,|\,F(x)\in \mathbb S^m_+\}
\end{equation}
where $f\colon\mathbb X\to\R$ and $F\colon\mathbb X\to\mathbb S^m$ are twice
continuously differentiable mappings and $\mathbb X$ is some Euclidean space.
Let $\mathcal F\subset\mathbb X$ be the feasible set of \eqref{eq:optimization_problem}.
For $\alpha\geq 0$, we introduce the generalized Lagrangian function
$\mathcal L^\alpha\colon\mathbb X\times\mathbb S^m\to\R$
associated with \eqref{eq:optimization_problem} by means of
\[
	\forall x\in\mathbb X\,\forall Y^*\in\mathbb S^m\colon\quad
	\mathcal L^\alpha(x,Y^*):=\alpha f(x)+\innerprod{Y^*}{F(x)}.
\] 
Furthermore, for $x\in\mathcal F$, we exploit the critical cone associated with
\eqref{eq:optimization_problem} given by
\[
	\mathcal C(x)
	:=
	\left\{
		u\in\mathbb X\,\middle|\,
		f'(x)u\leq 0,\,F'(x)u\in\mathcal T_{\mathbb S^m_+}(F(x))
	\right\}.
\]
Note that, due to \eqref{eq:tangent_cone_semidefinite}, this cone can be
computed explicitly as soon as an ordered eigenvalue decomposition of
$F(x)$ is at hand.
For $u\in\mathcal C(x)$ and $\alpha\geq 0$, 
the associated directional Lagrange multiplier set is given by
\[
	\Lambda^\alpha(x,u)
	:=
	\left\{Y^*\in \mathcal N_{\mathbb S^m_+}(F(x))\cap \{F'(x)u\}^\perp\,\middle|\,
		(\mathcal L^\alpha)'_x(x,Y^*)=0
	\right\},
\]
and this set can be computed via \eqref{eq:normal_cone_semidefinite}.

The following second-order sufficient optimality condition for
\eqref{eq:optimization_problem} can be distilled from the
more general result \cite[Theorem~3.3]{BeGfrYeZhouZhang1}
which has been proven via a straight contradiction argument,
and a direct proof of it, which is merely based on calculus rules for the
second subderivative, is stated in \cite[Theorem~5.2]{BenkoMehlitz2022c}.
A slightly less general result, which clearly
motivated the authors of \cite{BeGfrYeZhouZhang1}, can be found in
\cite[Theorem~7.1]{MohammadiMordukhovichSarabi2021}.
\begin{theorem}\label{thm:SOSC_general}
	Let $\bar x\in\mathcal F$ be chosen such that
	for each $u\in\mathcal C(\bar x)\setminus\{0\}$, there are
	$\alpha\geq 0$ and $Y^*\in\Lambda^\alpha(\bar x,u)$ such that
	\begin{equation}\label{eq:SOSC_estimate}
		(\mathcal L^\alpha)''_{xx}(\bar x,Y^*)[u,u]
		+
		\mathrm d^2\delta_{\mathbb S^m_+}(F(\bar x),Y^*)(F'(\bar x)u)
		>
		0.
	\end{equation}
	Then $\bar x$ is an essential local minimizer of second order
	for \eqref{eq:optimization_problem}, 
	i.e., there are $\varepsilon>0$ and $\beta>0$ such that
	\begin{equation}\label{eq:essential_local_minimizer}
		\forall x\in\mathbb B_\varepsilon(\bar x)\colon\quad
		\max\left(f(x)-f(\bar x),\dist_{\mathbb S^m_+}(F(x))\right)
		\geq
		\beta\norm{x-\bar x}^2.
	\end{equation}
	Particularly, $\bar x$ is a strict local minimizer of \eqref{eq:optimization_problem}.
\end{theorem}

It is clear by definition of the second subderivative that \eqref{eq:SOSC_estimate} can
only hold for some $u\in\mathcal C(\bar x)\setminus\{0\}$, $\alpha\geq 0$, and
$Y^*\in\Lambda^\alpha(\bar x,u)$ if $(\alpha,Y^*)\neq(0,O)$, i.e., the non-triviality
of the appearing generalized Lagrange multipliers is inherent.

We also note that the growth condition
\eqref{eq:essential_local_minimizer} is slightly more restrictive than 
\[
	\forall x\in\mathcal F\cap\mathbb B_\varepsilon(\bar x)\colon\quad
	f(x)-f(\bar x)
	\geq
	\beta\norm{x-\bar x}^2
\]
which is referred to as the second-order growth condition associated with
\eqref{eq:optimization_problem} at $\bar x$ in the literature.

In order to turn \eqref{eq:SOSC_estimate} into a valuable second-order optimality
condition, the appearing second subderivative of $\delta_{\mathbb S^m_+}$ 
has to be evaluated or, at least,
estimated from below. Exemplary, this strategy has been used in \cite{BenkoMehlitz2022c}
in order to infer second-order sufficient conditions in nonlinear second-order cone programming
and turned out to be much simpler than the more technical verification strategies from
\cite{BonnansRamirez2005,HangMordukhovichSarabi2020}.
Here, we present a similar analysis for nonlinear semidefinite programs.
As already remarked in \cite{BenkoMehlitz2022c}, 
obtaining second-order necessary optimality conditions
based on second subderivatives is often not reasonable since
this would come along with comparatively strong regularity conditions which are necessary
in order to get the calculus rules for second subderivatives working.

In the subsequent lemma, an explicit formula for the second subderivative of 
$\delta_{\mathbb S^m_+}$ is presented.
\begin{lemma}\label{lem:second_subderivative_semidefinite_cone}
	For each $Y\in\mathbb S^m_+$, $V\in\mathcal T_{\mathbb S^m_+}(Y)$, and
	$Y^*\in\mathcal N_{\mathbb S^m_+}(Y)\cap\{V\}^\perp$, we have
	\[
		\mathrm d^2\delta_{\mathbb S^m_+}(Y,Y^*)(V)
		=
		-2\ninnerprod{Y^*}{VY^\dagger V}.
	\]
\end{lemma}
\begin{proof}
	Let $Y=P^\top MP$ be an ordered eigenvalue decomposition of $Y$
	with orthogonal matrix $P\in\R^{m\times m}$ and diagonal matrix $M\in\mathbb S^m$
	as well as the index sets $\pi$ and $\omega$ as defined in \cref{sec:matrix_analysis}.
	From $Y^*\in\mathcal N_{\mathbb S^m_+}(Y)$, we find 
	$(Y^*)_{\pi\pi}^P=O$, $(Y^*)_{\pi\omega}^P=O$, 
	and $(Y^*)_{\omega\omega}^P\in\mathbb S^{|\omega|}_-$.
	Furthermore, $V\in\mathcal T_{\mathbb S^m_+}(Y)$ gives
	$V^P_{\omega\omega}\in\mathbb S^{|\omega|}_+$.
	From $\innerprod{Y^*}{V}=0$	and orthogonality of $P$, we have
	\[
		0
		=
		\innerprod{Y^*}{V}
		=
		\ninnerprod{(Y^*)^P}{V^P}
		=
		\ninnerprod{(Y^*)^P_{\omega\omega}}{V^P_{\omega\omega}}
	\]
	which gives $\innerprod{(Y^*)^P_{\omega\omega}}{V^P_{\omega\omega}}=0$.
	
	For given $V'\in\mathbb S^m$ and sufficiently small $t>0$,
	$M_{\pi\pi}+t(V')^P_{\pi\pi}$ is positive definite, and
	since $Y+tV'\in\mathbb S^m_+$ and $M+t(V')^P\in\mathbb S^m_+$ are
	equivalent by orthogonality of $P$, 
	\cref{lem:definiteness_block_matrix} can be used
	to infer that, for small enough $t>0$, $Y+tV'\in\mathbb S^m_+$
	equals
	\begin{equation*}
		(V')^P_{\omega\omega}
		-
		t(V')^P_{\omega\pi}\bigl[M_{\pi\pi}+t(V')^P_{\pi\pi}\bigr]^{-1}(V')^P_{\pi\omega}
		\in 
		\mathbb S^{|\omega|}_+.
	\end{equation*}
	Thus, from $(Y^*)^P_{\omega\omega}\in\mathbb S^{|\omega|}_-$, we find
	\begin{align*}
		\mathrm d^2\delta_{\mathbb S^m_+}(Y,Y^*)(V)
		&=
		\liminf\limits_{\substack{t\downarrow 0,\,V'\to V\\Y+tV'\in\mathbb S^m_+}}
		-\frac{2\innerprod{Y^*}{V'}}{t}
		=
		\liminf\limits_{\substack{t\downarrow 0,\,V'\to V\\Y+tV'\in\mathbb S^m_+}}
		-\frac{2}{t}\ninnerprod{(Y^*)^P_{\omega\omega}}{(V')^P_{\omega\omega}}
		\\
		&
		\geq
		\liminf\limits_{\substack{t\downarrow 0,\,V'\to V\\Y+tV'\in\mathbb S^m_+}}
		-2\innerprod{(Y^*)^P_{\omega\omega}}
		{(V')^P_{\omega\pi}\bigl[M_{\pi\pi}+t(V')^P_{\pi\pi}\bigr]^{-1}(V')^P_{\pi\omega}}
		\\
		&=
		-2\ninnerprod{(Y^*)^P_{\omega\omega}}{V^P_{\omega\pi}M_{\pi\pi}^{-1}V^P_{\pi\omega}}
		=
		-2\ninnerprod{(Y^*)^P}{V^PM^\dagger V^P}
		\\
		&
		=
		-2\ninnerprod{Y^*}{VY^\dagger V}.
	\end{align*}
	
	Finally, we construct particular sequences $\{t_k\}_{k\in\N}\subset(0,\infty)$
	and $\{V_k\}_{k\in\N}\subset\mathbb S^m$ which show that this lower estimate
	is sharp. Therefore, let $\{t_k\}_{k\in\N}\subset(0,\infty)$ be a null sequence
	such that $M_{\pi\pi}+t_kV^P_{\pi\pi}$ is invertible for each $k\in\N$.
	Define
	\[
		\Delta_k
		:=
		t_kV^P_{\omega\pi}\bigl[M_{\pi\pi}+t_kV^P_{\pi\pi}\bigr]^{-1} V^P_{\pi\omega}
	\]
	and
	\[
		V_k
		:=
		P^\top 
			\begin{bmatrix}
				V^P_{\pi\pi}&V^P_{\pi\omega}\\
				V^P_{\omega\pi}&V^P_{\omega\omega}+\Delta_k
			\end{bmatrix}
		P
	\]
	for each $k\in\N$. Clearly, we have $\Delta_k\to O$ which gives $V_k\to V$.
	By construction, we also have 
	\begin{align*}
		(V_k)^P_{\omega\omega}
		&=
		V^P_{\omega\omega}+\Delta_k
		\\
		&=
		V^P_{\omega\omega}
		+
		t_kV^P_{\omega\pi}\bigl[M_{\pi\pi}+t_kV^P_{\pi\pi}\bigr]^{-1} V^P_{\pi\omega}
		\\
		&=
		V^P_{\omega\omega}
		+
		t_k(V_k)^P_{\omega\pi}\bigl[M_{\pi\pi}+t_k(V_k)^P_{\pi\pi}\bigr]^{-1}(V_k)^P_{\pi\omega},
	\end{align*}
	and rearrangements lead to 
	\[
		(V_k)^P_{\omega\omega}
		-
		t_k(V_k)^P_{\omega\pi}\bigl[M_{\pi\pi}+t_k(V_k)^P_{\pi\pi}\bigr]^{-1}(V_k)^P_{\pi\omega}
		=
		V^P_{\omega\omega}\in\mathbb S^{|\omega|}_+.
	\]
	Thus, \cref{lem:definiteness_block_matrix} gives $Y+t_kV_k\in\mathbb S^m_+$ for each $k\in\N$.
	Reprising the above steps for the estimation of the lower limit and
	recalling $\innerprod{(Y^*)^P_{\omega\omega}}{V^P_{\omega\omega}}=0$, we find
	\begin{align*}
		\liminf\limits_{k\to\infty}
		-\frac{2\innerprod{Y^*}{V_k}}{t_k}
		&=
		\liminf\limits_{k\to\infty}
		-\frac{2}{t_k}\ninnerprod{(Y^*)^P_{\omega\omega}}{(V_k)^P_{\omega\omega}}
		\\
		&
		=
		\liminf\limits_{k\to\infty}
		-\frac{2}{t_k}\innerprod{(Y^*)^P_{\omega\omega}}
		{V^P_{\omega\omega}+t_kV^P_{\omega\pi}\bigl[M_{\pi\pi}+t_kV^P_{\pi\pi}\bigr]^{-1} V^P_{\pi\omega}}
		\\
		&=
		\liminf\limits_{k\to\infty}
		-2\innerprod{(Y^*)^P_{\omega\omega}}{V^P_{\omega\pi}\bigl[M_{\pi\pi}+t_kV^P_{\pi\pi}\bigr]^{-1} V^P_{\pi\omega}}
		\\
		&=
		-2\ninnerprod{Y^*}{VY^\dagger V}.
	\end{align*}
	This already completes the proof.
\end{proof}

Let us note that the assertion of \cref{lem:second_subderivative_semidefinite_cone}
has been proven in \cite[Example~3.7]{MohamadiSarabi2020} with the aid of some
deeper results from \cite{Torki1999} addressing variational properties of
eigenvalue functions. In contrast, our proof is rather elementary.

Combining this result with \cref{thm:SOSC_general}, we obtain fully explicit
second-order sufficient optimality conditions for \eqref{eq:optimization_problem}.

\begin{corollary}\label{thm:SOSC_semidefinite}
	Let $\bar x\in\mathcal F$ be chosen such that
	for each $u\in\mathcal C(\bar x)\setminus\{0\}$, there are
	$\alpha\geq 0$ and $Y^*\in\Lambda^\alpha(\bar x,u)$ such that
	\begin{equation}\label{eq:SOSC_semidefinite}
		(\mathcal L^\alpha)''_{xx}(\bar x,Y^*)[u,u]
		>
		2\ninnerprod{Y^*}{(F'(\bar x)u)F(\bar x)^\dagger (F'(\bar x)u)}.
	\end{equation}
	Then $\bar x$ is an essential local minimizer of second-order for the associated
	optimization problem \eqref{eq:optimization_problem}.
\end{corollary}

Let us point the reader's attention to the simplicity of the above arguments which have 
been used to obtain this second-order optimality condition. \cref{thm:SOSC_general} is
proven via a standard contradiction argument. Further, the computation of the
appearing second subderivative of $\delta_{\mathbb S^m_+}$ is completely elementary
and relies on the standard approach of working with an ordered eigenvalue decomposition.
In \cite[Theorem~9]{Shapiro1997} and \cite[Section~5.3.5]{BonSh00}, related 
second-order sufficient conditions, based on the same expression for the curvature
term, i.e., the right-hand side in \eqref{eq:SOSC_semidefinite}, 
but with a weaker growth condition were obtained using the 
theory of second-order tangent sets. This approach is much more technical and relies
on deeper mathematics such as second-order directional differentiability of the
smallest eigenvalue of a positive semidefinite matrix.

\section{Concluding remarks}\label{sec:conclusions}

In this note, we computed the second subderivative of the indicator function associated
with the cone of all positive semidefinite matrices, and
this finding was used to obtain second-order sufficient optimality conditions in
nonlinear semidefinite optimization. 
This procedure recovered the findings from \cite{Shapiro1997} in elementary way. 
In the future, it needs to be studied whether this second-order sufficient condition
can be employed beneficially in numerical optimization like in \cite{HangMordukhovichSarabi2022}
where local analysis of a multiplier-penalty method associated with second-order cone programs
is investigated.
Furthermore, it seems reasonable to check whether our approach to second-order sufficient
conditions yields comprehensive results when applied to optimization problems with
semidefinite cone complementarity constraints,
see e.g.\ \cite{DingSunYe2014,LiuPan2022,WuZhangZhang2014}.
Finally, we note that
\[
	\mathbb S^m_+
	=
	\{
		Y\in\mathbb S^m\,|\,
		\forall v\in\R^m\colon\,v^\top Yv\geq 0
	\}
\]
holds, so $\mathbb S^m_+$ is a special instance of the closed, convex cone
\[
	\mathbb S^m_+(K)
	:=
	\{
		Y\in\mathbb S^m\,|\,
		\forall v\in K\colon\,v^\top Yv\geq 0
	\}
\]
where $K\subset\R^m$ is an arbitrary closed, convex cone.
In the literature, $\mathbb  S^m_+(K)$ is referred to as the set-semidefinite or
set-copositive cone
associated with $K$, and for $K:=\R^m_+$, the popular copositive cone is obtained,
see \cite{Bomze2012,Burer2015,Duer2010,DuerRendl2021} for further information about this
cone and applications of copositive optimization.
Following the approach of this note, it might be possible to obtain second-order sufficient
conditions for nonlinear optimization problems involving $\mathbb S^m_+(K)$. 
However, it is well known that the variational geometry of 
$\mathbb S^m_+(K)$ is much more challenging for general $K$ than for $K:=\R^m$, so the
necessary computations might be much more involved than the ones 
from \cref{lem:second_subderivative_semidefinite_cone}.


%

\end{document}